\documentclass[reqno]{amsart}

\usepackage{graphics,color}

\usepackage{amssymb,amsmath,color}
\usepackage{amssymb,amsmath}
\theoremstyle{definition}

\newtheorem{theorem}{Theorem}[section]
\newtheorem{lemma}[theorem]{Lemma}

\newtheorem{cor}[theorem]{Corollary}

\newtheorem{prop}[theorem]{Proposition}
\def\l{{\lambda}}

\newcommand{\BQ}{{\mathbb Q}}
\newcommand{\BZ}{{\mathbb Z}}

\newcommand{\al}{\alpha}
\newcommand{\G}{{\Gamma}}
\newcommand{\SL}{{SL_2(\BZ)}}

      \let\d\delta
  
  \let\l\lambda

\def\C{\mathbb C}
\def\G{\mathbf G}
\def\A{\mathbf A}

\def\E{\mathcal E}

\def\CF{\mathcal F}

\def\CP{\mathcal P}

\def\Tr{{\rm Tr}}
\def\Gal{{\rm Gal}}
\def\Frob{{\rm Frob}}
\def\F{{\mathbb F}}
  
\def\Z{{\mathbb Z}}

\def\Q{{\mathbb Q}}

\def\G{\Gamma}

\DeclareRobustCommand{\qed}{%
  \ifmmode % if math mode, assume display: omit penalty etc.
  \else \leavevmode\unskip\penalty9999 \hbox{}\nobreak\hfill
  \fi
  \quad\hbox{$\square$} \\}

\begin{document}

\title{ON ATKIN and SWINNERTON-DYER CONGRUENCE RELATIONS (2)}

\author{A.O.L. Atkin}
\address{Department of Mathematics \\ University of Illinois at Chicago\\
Chicago, IL 60637\\USA} \email{aolatkin@math.uic.edu}

\author{Wen-Ching Winnie Li}
\address{Department of Mathematics \\Pennsylvania State University\\
University Park, PA 16802\\USA} \email{wli@math.psu.edu}

\author{Ling Long}
\address{Department of Mathematics\\Iowa State University\\ Ames, IA 50011 \\USA}
\email{linglong@iastate.edu}

\thanks{The research of the second author was supported in part by an NSA
grant \#MDA904-03-1-0069 and an NSF grant \#DMS-0457574. Part of the
research was done when she was visiting the National Center for
Theoretical Sciences in Hsinchu, Taiwan. She would like to thank the
Center for its support and hospitality. The third author was
supported in part by an NSF-AWM mentoring travel grant for women.
She would further thank the Pennsylvania State University and the
Institut des Hautes \'Etudes Scientifiques for their hospitality.}

%\date{April 14, 2006}

\maketitle

\begin{abstract}
In this paper we give an example of a noncongruence subgroup whose
three-dimensional space of cusp forms of weight 3 has  the following
properties. For each of the four residue classes of odd primes
modulo 8 there is a basis whose Fourier coefficients at infinity
satisfy a three-term Atkin and Swinnerton-Dyer congruence relation,
which is the $p$-adic analogue of the three-term recursion satisfied
by the coefficients of classical Hecke eigenforms. We also show that
there is an automorphic $L$-function over $\mathbb Q$ whose local
factors agree with those of the $l$-adic Scholl representations
attached to the space of noncongruence cusp forms.

%\subclass{11F11, 11F33}

\end{abstract}

\section{Introduction}
The cusp forms of weight $k$ for congruence subgroups of the
classical modular  group $SL_2(\BZ)$ form a vector space of finite
dimension, and have a basis of  forms which are simultaneous
eigenfunctions of almost all the Hecke operators. In terms of the
series expansions at the cusp infinity, we can write such a form
as $\sum _{n\ge 1} a(n) w^n$, where $w=e^{2 \pi i z/\mu}$ is the
local uniformizer, and the Fourier coefficients $a(n)$ satisfy
\begin{equation}\label{eq:1}
     a(np) - A(p)a(n) + B(p) a(n/p) =0
\end{equation}
for all $n\ge 1$ and almost all primes $p$; as usual the
number-theoretic function $a(x)$ is defined to be zero if $x$ is not
a rational integer. The values of $a(n)$, $A(p)$, and $B(p)$ lie in
an algebraic number field, and
\begin{equation*}
                |A(p)|\le 2p^{(k-1)/2}
\end{equation*} and
\begin{equation}\label{eq:2}
                |B(p)|= p^{k-1}
\end{equation}
 for almost all $p$. If we normalize our basis by setting $a(1)=1$, it is
clear that we have
 $A(p)=a(p)$ for all relevant $p$; there is also a Dirichlet character
 $\chi$ such that
 $B(p) = \chi(p)p^{k-1}$.

It has been known for a long time that noncongruence subgroups of
the  modular group exist, and the finite-dimensional vector space of
their cusp forms of weight $k$ will still have a basis of forms with
series expansions and Fourier coefficients $a(n)$. However the
classical theory of the Hecke operators collapses so that no
three-term relation like \eqref{eq:1} above can be expected to hold
identically.

The first systematic investigation of noncongruence forms was
carried out  by Atkin and Swinnerton-Dyer \cite{a-sd}. First they
evolved an effective technique for constructing forms and functions
on noncongruence subgroups. Then for a small number of explicit
examples they found an appropriate analogy of \eqref{eq:1}, without
making any formal conjectures, and without proof except for the
isolated case of weight 2 and dimension 1.

For almost all $p$ they replaced the equality in \eqref{eq:1} by
congruence modulo a suitable power of $p$ (which we call a 3-term
ASD congruence here), with a $p$-adic basis of forms (so that the
$a(n)$ were $p$-adic); the $A(p)$ were however apparently algebraic
in general, although in different number fields for each $p$ (and of
course could no longer be identified as $a(p)$), and the $B(p)$ were
$p^{k-1}$ times a root of unity. Using $p$-adic Lie groups, Cartier
\cite{car71} established the 3-term ASD congruence for noncongruence
weight-2 forms.

A major advance was made by Scholl some fifteen years later
\cite{sch85b}. To the $d$-dimensional space of cusp forms of weight
$k (\ge 3)$ of a noncongruence subgroup whose modular curve was
defined over $\BQ$, Scholl associated a family of $2d$-dimensional
$l$-adic representations of the Galois group over $\BQ$ such that
the characteristic polynomials $P_p(T)$ of the Frobenius elements at
almost all $p$ were degree $2d$ polynomials over $\BZ$, independent
of $l$, and all the cusp forms in this space with algebraic Fourier
coefficients satisfied a $(2d+1)$-term congruence relation given by
$P_p(T)$. Moreover, all zeros of $P_p(T)$ had absolute value equal
to $p^{(k-1)/2}$.  In particular, he showed that the 3-term ASD
congruence was valid for the case of a 1-dimensional space. When $d
> 1$, while Scholl did not refine the $(2d+1)$-term congruence
to establish the 3-term ASD congruences, he did point out in
Sections 5.6 and 5.8 of \cite{sch85b} that if the Frobenius
matrices of the Scholl's representation are diagonalizable (which
are expected conjecturally) and the primes are mostly ordinary
(which are expected for forms with no CM), then the 3-term ASD
congruences will follow.

A related problem is the interpretation of the algebraic numbers
$A(p)$ and $B(p)$.
 In the case $d=1$, these are of course rational integers, and Scholl
 established
a number of examples (in \cite{sch88} and \cite{sch04}) where $A(p)$
were the Fourier coefficients of a congruence modular form of weight
$k$ and $B(p) = p^{k-1}$. This also holds for weight 3 noncongruence
cusp forms associated to K3 surfaces over $\Q$, as shown in
\cite{lly05}.

More recently Li, Long, and Yang \cite{lly05} have given an example
with dimension $2$ and weight 3, where the eigenforms are the same
for all $p$ (so that the $a(n)$ in \eqref{eq:1} are in a number
field and not merely $p$-adic), and the numbers $A(p)$ are the
Fourier coefficients of two congruence weight-3 forms of character
$\chi$, and $B(p) = \chi(p) p^2$. This is as much as could possibly
be hoped for, and is clearly exceptional. Another example in the
same vein is obtained by Fang et al in \cite{F5}.

In the present paper we give an example with dimension $3$ which can
be decomposed into a 1-dimensional space and a 2-dimensional space
where the basis of the  1-dimensional space satisfies the 3-term ASD
congruence with the $A(p)$ being coefficients of an explicit
weight-3 cusp form. In the 2-dimensional space,
\begin{itemize}
\item[1)] for each of the four residue classes of odd primes
modulo 8 there is a basis (again over a number field and not merely
$p$-adic) of two forms, for each of which the $a(n)$ enjoy the
3-term ASD congruence;
\smallskip

\item[2)] the corresponding $A(p)$ are again the Fourier
coefficients of congruence modular forms, and $B(p) = \pm p^2$.
\end{itemize} Our main result is recorded in Theorem \ref{ASDexplicit}.
While not as simple as the example in \cite{lly05}, it still
exhibits far more structure than the general situation. It may be
significant that both these examples arise from noncongruence
groups which are normal subgroups of congruence subgroups of genus
zero, with a cyclic factor group.

We shall show that there is an automorphic $L$-function over $\BQ$
whose local factors agree with those of the $l$-adic Scholl
representations attached to the space of weight-3 cusp forms of
our noncongruence group. This exemplifies the Langlands
philosophy, which predicts that $l$-adic representations that are
motivic and ramified at finitely many places should relate to
representations of $GL_n$. When this happens, we say that the
$l$-adic representations are modular. It is worth pointing out
that the modularity is proved by first restricting both Scholl's
$l$-adic representation and the $l$-adic representation attached
to the automorphic form to the Galois group of $\mathbb Q(i)$ and
then comparing them using Livn\'e's criterion \cite{liv87}, since
both residual representations have even trace. This is different
from the method in \cite{lly05}, where the Galois group was over
$\mathbb Q$ and Serre's criterion \cite{ser_1} was used. To the
best of our knowledge, the modularity technique introduced by
Wiles is not yet applicable when the base field is $\Q(i)$.
\smallskip

This paper is organized as follows. In Sec. 2 we introduce the
 normal subgroups $\G_n$ of $\G^1(5)$ and study
 properties of the operators involved.
 Using an explicit model of an elliptic surface over the modular
 curve of interest to us, we construct in Sec. 3 a family of
$l$-adic representations of the Galois group over $\mathbb Q$,
which are isomorphic to Scholl representations. Taking advantage
of the concrete model, we compute the traces of the Frobenius
elements, and obtain the determinants of the representations. The
case of $\G_2$, studied in \cite{lly05} using geometry, is
reviewed in Sec. 4.

The rest of the paper concerns the group $\G_4$, whose
$3$-dimensional space of weight-3 cusp forms is to be analyzed. As
shown by Scholl \cite{sch85b}, there is a 7-term ASD congruence at
$p$, obtained by comparing the $l$-adic theory and the $p$-adic
theory, more precisely, from the agreement of the characteristic
polynomials of the Frobenius element at $p$ ~in both theories. In
order to obtain 3-term ASD congruences, we need to decompose the
vector space of each side into three 2-dimensional invariant
subspaces and prove the agreement of the characteristic
polynomials on the 2-dimensional subspaces. This is achieved by
taking eigenspaces of suitably chosen operators which commute with
the action of the Frobenius at $p$. These operators, which play
the role of the Hecke operators on congruence forms, are
constructed using the symmetries of the elliptic surface.

The details are carried out as follows. In Sec. 5 we decompose
each side into the direct sum of a 2-dimensional $(+)$-space and a
4-dimensional $(-)$-space. On the $l$-adic side, this gives the
decomposition of the Scholl representation $\rho_l$ into the sum
of $\rho_{l,+}$ and $\rho_{l,-}$. The $(+)$-spaces can be easily
identified as arising from the group $\G_2$, and hence the
agreement of the characteristic polynomials. The operators used to
further decompose the $(-)$-spaces into the sum of two
2-dimensional subspaces are introduced in Sec. 6. On the
$(-)$-space of the $l$-adic side, these operators yield, for each
odd prime $p \ne l$, a factorization of the characteristic
polynomial of the Frobenius element
 at $p$ as a product of two degree-2
polynomials of a specific type. Aided by Magma, we are able to
obtain these factors explicitly for $p = 3, 7,$ and $ 13$. By
comparing the actions of these operators on both $(-)$-spaces, we
show in Sec. 7 that there is a basis of the $(-)$-space satisfying
the 3-term ASD congruence at $p$ given by the two degree-2 factors
of the characteristic polynomial at $p$ obtained in Sec. 6.
Finally in Sec. 8 we identify a cusp form $f'$ of weight 3 and an
idele class character $\chi$ of $\mathbb Q(i)$ of order 4, which
yields a cusp form $h(\chi)$ for $GL_2(\mathbb Q)$, and prove that
the $l$-adic representation attached to $f' \times h(\chi)$ and
$\rho_{l,-}$ have the same semi-simplification. This establishes
the modularity of $\rho_{l.-}$ and gives explicit ASD congruences
at the same time.

\section{The groups and operators}
\subsection{The groups $\G_n$} Listed below are the cusps of
$\pm\Gamma^1(5)$ and a choice of generators of their stabilizers in
$SL_2(\Z)$:
\begin{center}
\begin{tabular}{c c c c c c}
cusps of $\pm\Gamma^1(5)$ \qquad&cusp widths& \qquad & generators& of& stabilizers\\
\hline $\infty$ \qquad &5&\qquad
 & $\zeta$&=&$\displaystyle{\left(\begin{matrix}1&5\\0&1\end{matrix}
                   \right)}$\\
$0$ \qquad&1&\qquad
 & $\theta$&=&$\displaystyle{\left(\begin{matrix}1&0\\-1&1\end{matrix}
                 \right)}$\\
$-2$ \qquad&5&\qquad
 & $A\zeta
A^{-1}$&=&$\displaystyle{\left(\begin{matrix}11&20\\-5&-9\end{matrix}
                 \right)}$\\
$-\frac52$ \qquad&1&\qquad
 & $A\theta
A^{-1}$&=&$\displaystyle{\left(\begin{matrix}11&25\\-4&-9\end{matrix}
                 \right)}$\\
\end{tabular}
\end{center}

\noindent Here $A = \begin{pmatrix} -2&-5\\1&2
\end{pmatrix}$. Therefore the group
$\Gamma^1(5)$ is generated by $\zeta$, $\theta$, $A\zeta A^{-1}$,
$A\theta A^{-1}$ with the relation
$$(A\theta
A^{-1})(A\zeta A^{-1})\theta \zeta = I_2.$$ Let $E_1(z)$ and
$E_2(z)$ be two weight-3 Eisenstein series of $\G^1(5)$ with
integral coefficients, vanishing at all cusps except at $\infty$ and
$-2$, respectively; these series are constructed explicitly in
\cite{lly05}. The function $t=\frac{E_1(z)}{E_2(z)}$ generates the
field of meromorphic functions on $X^1(5)$, the modular curve
associated to $\G^1(5)$.

 For $n \ge 2$, let $\varphi_n$ denote a homomorphism  from $\G^1(5)$ to $\C^{\times}$
which sends $\zeta, A\zeta A^{-1}$, to $\omega_n, \omega_n^{-1}$,
respectively,  and  the other generators to $1$, where
$\omega_n=e^{2\pi i/n}$. Denote the kernel of
 $\varphi_n$ by $\G_n$; it is an index-$n$ subgroup of $\G^1(5)$, generated by
$$\zeta^n, A\zeta^n A^{-1}, \zeta^i \theta \zeta^{-i},
\zeta^{i} A \theta A^{-1} \zeta^{-i}, \text{ for} \quad
i=0,\cdots,n-1,$$

\noindent subject to the relation
\begin{eqnarray*}
  (A\theta A^{-1})(A\zeta^n A^{-1}) \theta (\zeta A\theta
A^{-1}\zeta^{-1}) (\zeta\theta\zeta^{-1})\cdots\\
(\zeta^{n-1} A\theta A^{-1}\zeta^{-(n-1)})
(\zeta^{n-1}\theta\zeta^{-(n-1)})\zeta^n = I_2.
\end{eqnarray*}

The modular curve of $\G_n$, denoted by $X_n$, is of genus 0. It is
an $n$-fold cover of the modular curve of $\G^1(5)$, ramified only
at the cusps $\infty$ and $-2$, both with ramification degree $n$. A
generator for the field of meromorphic functions of $X_n$ is
\begin{equation}
  t_n=\sqrt[n]{\frac{E_1}{E_2}}.
\end{equation}

 A cusp form is called $n$-integral if its Fourier coefficients
are algebraic and integral outside the places dividing $n$.

\begin{prop}
 \begin{enumerate}
  \item  $\G_n$ is a normal subgroup of $\G^1(5)$ with the quotient $\G^1(5)/\G_n$
  cyclic of order $n$ generated by $\zeta$. The characters of this quotient
  are  ${\varphi_n}^j$, $1 \le j \le n$, where ${\varphi_n}$ is as above.
  \item $\dim S_3(\G_n)=n-1$.
  \item An $n$-integral basis of $S_3(\G_n)$ consists of $
  \sqrt[n]{E_1^{n-j}(z)E_2^{j}(z)}$, $1 \le j \le n-1$.
  \end{enumerate}
\end{prop}
\begin{proof}
    Since $\G_n$ is the kernel of the homomorphism $\varphi_n$, it is a
    normal subgroup of $\G^1(5)$. It follows from the definition of
   $\varphi_n$ that the quotient $\G^1(5)/\G_n$ is the cyclic group
   generated by $\zeta$, and hence the characters are as described in (1).

  The group $\G_n$ has no elliptic elements and has $2+2n$ cusps,
  of which two have width $5n$, and $2n$ have width 1. Moreover
  $-I\notin \G_n$. Using the
  dimension formula in \cite[Theorem2.25]{shim1},  one concludes $\dim
  S_3(\G_n)=n-1.$

  It is easy to verify that $\sqrt[n]{E_1^{n-j}(z)E_2^{j}(z)}$,
  $1 \le j \le n-1$, are all weight-3
  cusp forms on $\G_n$ and are linearly
  independent, as seen from their Fourier expansions, which we
  normalize as
  $$ \sqrt[n]{E_1^{n-j}(z)E_2^{j}(z)}=e^{(2\pi iz)j/5n}(1+ \sum_{r \ge 1} c_r e^{2\pi izr/5}),$$
  where the coefficients $c_r$ are $n$-integral, since $E_1,\, E_2\in \Z[[e^{2\pi iz/5}]].$
\end{proof}

\subsection{Operators} The matrices $A=\begin{pmatrix}   -2&-5\\1&2
\end{pmatrix}$ and $\zeta =\begin{pmatrix}1&5\\0&1
\end{pmatrix}$ normalize $\G^1(5)$ and all the $\G_n$. Note that $A$ is the diamond
operator $\langle 2 \rangle$. The elliptic modular surface
$\E_{\G^1(5)}$ over $\G^1(5)$ in the sense of Shioda \cite{shio1}
has a model given by
\begin{multline}\label{eq:E1(5)}
  y^2=t(x^3-\frac{1+12t+14t^2-12t^3+t^{4}}{48t^2}x
  \\+\frac{1+18t+75t^2+75t^{4}-18t^{5}+t^{6}}{864t^3} ),
\end{multline}
where $t$ is the parameter for the genus zero modular curve
$X^1(5)$. The action of $A$ on this elliptic surface is
$$A: \quad (x,y,t)\mapsto (-x,y/t,-1/t).$$ It is a morphism of
order 4 defined over $\Q$. The operator
$\zeta$ acts trivially on $\E_{\G^1(5)}$. The actions of $A$ and
$\zeta$ on the modular curve $X_n$ for $\G_n$ are described in
terms of
 the function $t_n$ as follows:
\begin{eqnarray}
  A(t_n)=\frac{\omega_{2n}}{t_n},& \quad &
  \zeta (t_n)=\omega_{2n}^{-2} t_n,
\end{eqnarray} where $\omega_{2n}= e^{2\pi i/2n}$ is a primitive $2n$-th root of unity.
For even $n$, an equation of the elliptic modular surface $\E_n$
attached to $\G_n$ in the sense of Shioda \cite{shio1} is
\begin{multline}\label{eq:E_n}
 y^2=x^3-\frac{1+12(t_n^{n}-t_n^{3n})+14t_n^{2n}+t_n^{4n}}{48t_n^{2n}}x\\
  +\frac{1+18(t_n^{n}-t_n^{5n})+75(t_n^{2n}+t_n^{4n})+t_n^{6n}}{864t_n^{3n}}.
 \end{multline}

\section{$l$-adic representations}
\subsection{Two $l$-adic representations}
 For any field $K$, write $G_{K}$ for $\Gal(\overline{K}/K)$.
%For a prime $l$,
 %let $W_{n,l}$  denote the $\Q_l$ vector space for Scholl's
%$l$-adic representation $$\rho_{n,l}: G_{\Q} \rightarrow
%Aut(W_{n,l})$$ attached to the space of weight $3$ cusp forms of
%$\G_n$.
As in \cite{lly05}, for a prime $l$, we introduce a family of
$l$-adic representations $\rho_{n,l}^*$ of $G_\Q$ attached to the
space of weight $3$ cusp forms of $\G_n$ by using  the explicit
model (\ref{eq:E_n}) of the elliptic surface $\E_n$ over $X_n$ as
follows.

Denote by $X_n^0$ the modular curve $X_n$ with the cusps removed.
Let $h: \E_n \rightarrow X_n$ be the natural elliptic fibration
endowed with $\E_n$. Let $h^0: \E_n \rightarrow X_n^0$ be its
restriction to $X_n^0$, which is a smooth map. For any prime $l$,
one obtains a sheaf
\[
\CF_l=R^1h^0_* \Q_l
\]
on $X_{n}^0$.  Here $\Q_l$ is the constant sheaf on the elliptic
surface $\E_n$ and $R^1$ is the derived functor. The inclusion map
$i: X_n^0 \rightarrow X_n$ then transports the sheaf on $X_n^0$ to
a sheaf $i_*\CF_l$ on $X_{n}$. The action of $G_{\Q}$ on the
$\Q_l$-space
\begin{equation}\label{e:4.2}
%W_l=H^1(X_{\G'}[1/Ml]\otimes\overline\Q,i_*\CF_l)
 W_{n,l}=H^1(X_{n}\otimes\overline\Q,i_*\CF_l)
\end{equation}
\noindent defines an $l$-adic representation, denoted by
$\rho_{n,l}^{*}$, of the Galois group $G_\Q$.

Scholl's representation $\rho_{n,l}$ of $G_{\Q}$
 attached to $S_3(\G_n)$ in \cite{sch85b}
is constructed by first choosing an auxiliary modular curve $X(N)$
and pulling the universal elliptic curve on $X(N)$ to the fibre
product of $X(N)$ with $X_n$, going through the same process as
above with $X_n$ replaced by its fibre product with $X(N)$, and at
the end taking the $SL(\mu_n\times \Z/N)$-invariant
 part of the cohomology space
to be the representation space of $\rho_{n,l}$ (cf. (2.1.2) of
\cite{sch85b}).
 Since
\eqref{eq:E1(5)} is an algebraic model for the universal elliptic
curve $\E_{\G^1(5)}$ of $\G^1(5)$ on $X_{\G^1(5)}$, we may use
$X_{\G^1(5)}$ as the auxiliary curve in Scholl's construction.
Then $\E_n$ is the pullback of $\E_{\G^1(5)}$ to the fibre product
of $X_n$ with $X_{\G^1(5)}$, which is isomorphic to $X_n$.
Therefore the two representations are isomorphic; we shall remove
$^*$ from now on.

 When $n=2$, it is shown in \cite{b-s} that the K3
surface \begin{equation}\label{eq:S}\mathfrak{S}:
t_2^2=uv(1-u)(v+1)(u-v)
\end{equation}  has
an elliptic fibration
\begin{equation}\label{eq: 8} y^2+(1-t_2^2)xy-t_2^2y=x^3-t_2^2x^2\end{equation}
over $X_2$. Furthermore,  the fibration defined by \eqref{eq: 8}
is birationally isomorphic (over $\Q$) to the fibration over
$\E_2$ defined by \eqref{eq:E_n}. The space $W_{2,l}$ constructed
above can be embedded into $H^2(\mathfrak{S},\Q_l)$, which shows
that the representation $\rho_{2,l}$ is unramified at 5 if $l \ne
5$. Since $\mathfrak{S}$ is obtained as a double cover of $\mathbb
P^2$ branched over six lines with branch locus shown in \cite[Fig.
1]{b-s}, we see that the zero and pole locus of $t_2=0$ is a
divisor on \eqref{eq:S} with normal crossings, which extends to a
divisor with relative normal crossing over $\Z_5$. When $n=4$,
$\E_4$ is obtained from $\E_2$ by replacing $t_2$ by $t_4^2$.
Therefore $\rho_{4,l}$ is also unramified at 5.

\subsection{Traces}\label{subsec:trace}

Let $q=p^r$ where $r$ is a positive integer. Denote by $\Frob_p$
the action of the Frobenius element at $p$ on an elliptic surface.

The  surface $\E_{\G^1(5)}$ is a
rational surface with vanishing $H^{2,0}(\E_{\G^1(5)},\Z)$ (i.e.,
$S_3(\G^1(5))= 0$). Hence $\dim W_{1,l}=0$. So for all
 odd prime powers $p^r$ we have
$$\Tr \rho_{1,l}(\Frob_{p^r})=0.$$

 We proceed to compute $\Tr
\rho_{n,l}(\Frob_{p^r})$ for $n = 2, 4$. The method is the same as
that used in Section 5 of \cite{lly05}. %Here we briefly review the
%key ingredients.
By the Lefschetz fixed point theorem,
$$\Tr \rho_{n,l}(\Frob_q) = - \sum_{t \in X_n(\F_q)} \Tr_n ~t,$$
\noindent where $\Tr_n ~t$ is the trace of the $\Frob_q$
restricted to the stalk at $t$ of $i_* \CF_l$. Similarly,
$$\Tr \rho_{1,l}(\Frob_q)= - \sum_{t \in X^1(5)(\F_q)} \Tr ~t.$$

 Recall that all the modular curves $X_n$
have genus zero, hence $X_n(\F_q)= \mathbb P^1(\F_q)$ can be
identified as $\F_q \cup \{\infty\}$. The curve $X_n$ is an
$n$-fold cover of $X^1(5)$ and so is $\E_n$ over $\E_{\G^1(5)}$
under the map $(x, y, t_n) \mapsto (x, y, t_n^n)$.  As such, the
fibre over $\infty$ of $X_n(\F_q)$ is the same as the fibre over $
\infty$ of $X^1(5)(\F_q)$ since $\E_{\G^1(5)}$ is semistable, and
hence has the same trace. Moreover, $\Tr_n t_n$ at $t_n \in
X_n(\F_q)$ is equal to $\Tr ~t_n^n$ at $t_n^n \in X^1(5)(\F_q)$.
The operator $A$ on $X^1(5)$ is $\Q$-rational; it induces a
bijection between the rational points on the fibre at $t \ne 0,
\infty$ and those at $-1/t$.

\begin{lemma}
  When $n=2,4$, for all prime powers $q=p^r$ with $q \equiv 3 \mod 4$,
  $$\Tr \rho_{n,l}(\Frob_q)=0.$$
\end{lemma}
\begin{proof} Note that %in this case
$-1$ is not a quadratic residue in $\F_q$. So $A$ induces a
bijection between the $\F_q$-rational points on the fibres at the
points of $X^1(5)(\F_q)$ parametrized by the quadratic residues of
$\F_q$ and those parameterized by the non-residues. Hence %Therefore

\begin{eqnarray*}
\Tr \rho_{2,l}(\Frob_q) & = - \sum_{t_2 \in X_2(\F_q)} \Tr_2 ~t_2
 = - \sum_{t \in X^1(5)(\F_q)} \Tr ~t^2 \\
 & = - \sum_{t \in
X^1(5)(\F_q)} \Tr ~t = \Tr \rho_{1,l}(\Frob_q)= 0.
\end{eqnarray*}
When $n=4$,
  note that $(\F_q^\times)^4 = (\F_q^\times)^2$ since $4 \nmid q-1$. So
\begin{eqnarray*}
  \Tr \rho_{4,l}(\Frob_q)=- \sum_{t_4 \in X_4(\F_q)} \Tr_4 ~t_4
  = - \sum_{t_2 \in X_2(\F_q)} \Tr_2 ~t_2 = \Tr \rho_{2,l}(\Frob_q)
  = 0.
  \end{eqnarray*}
\end{proof}

\begin{lemma}
 For all prime powers $q=p^r$ with $q \equiv 5 \mod 8$,
  $$\Tr \rho_{4,l}(\Frob_q)=\Tr \rho_{2,l}(\Frob_q).$$
  \end{lemma}
\begin{proof}
   The operator $A$ induces a
   bijection between $\F_q$-rational points on the fibres at
  the points of $X^1(5)(\F_q)$ parametrized by $(\F_q^\times)^4$
and those parameterized by $-(\F_q^\times)^4$. Observe that
$(\F_q^\times)^2 = (\F_q^\times)^4 \cup -(\F_q^\times)^4$ when
 $q \equiv 5 \mod 8$ since $-1$ is a quadratic residue, but not
 a 4-th power residue. Therefore
 \begin{eqnarray*}
\Tr \rho_{4,l}(\Frob_q) = - \sum_{t \in X^1(5)(\F_q)} \Tr ~t^4
  = - \sum_{t \in X^1(5)(\F_q)} \Tr ~t^2 = \Tr
  \rho_{2,l}(\Frob_q).
  \end{eqnarray*}
%\noindent as desired.
\end{proof}

\begin{lemma}\label{lem:ap=0mod8ifp=1mod8}
  For all prime powers $q=p^r$ with $q \equiv 1 \mod 8$, $$\Tr \rho_{4,l}(\Frob_q)-\Tr
  \rho_{2,l}(\Frob_q) \equiv 0 \mod 4.$$
\end{lemma}
\begin{proof}In this case, the field $\F_q$ contains 8th roots of unity.
Denote by $i$ an element in $\F_q$ of order 4.
  The desired statement is equivalent to
  $$\sum_{t\in X^1(5)(\F_q), ~t\neq 0,\infty} (\Tr \, t^2-\Tr \,t^4) \equiv
  \sum_{t\in X^1(5)(\F_q), ~t\neq 0,\infty} (\Tr\, t^2) \equiv 0 \mod
  4.$$ Owing to the symmetry given by $A$, we have
  $ \Tr \, t^2=\Tr \, (-t)^2 = \Tr (-1/t^2)$. This implies that
    four distinct values of $t \in \F_q^\times$  give rise to the same
    trace, except when $t^2 = -1/t^2$, that is, $t^4 = -1$. Thus the
    statement is reduced to
    $$\sum_{t \in X^1(5)(\F_q), t^4 = -1} \Tr\, t^2 \equiv 0 \mod
  4,$$ which is equivalent to $\Tr\, i - \Tr (-i)$ being even. In terms of the
  model for the elliptic surface
  over $X^1(5)$, this means checking the difference of the
  number of $\F_q$-rational points on $Y^2 = X^3 + \frac{1+2i}{4}X$
  and $Y^2 = X^3 + \frac{1-2i}{4}X$, which is obviously even.
\end{proof}

By running a Magma program calculating the traces of $\Frob_{p^r}$
on each fibre for small values of $r$ and $p=3,7,13, 17$, we
obtain the following characteristic polynomials, independent of
\textbf{$l \ne p$}~:
\begin{center}

\begin{tabular}{|c|c| c|}
 \hline
  $p$&\text{char. poly. for } $\rho_{2,l}(\Frob_p)$&\text{char. poly. for } $\rho_{4,l}(\Frob_p)$  \\
  \hline
  3&$x^2-3^2$&$(x^2-3^2)(x^4-10x^2+3^4)$\\
  \hline
  7&$x^2-7^2$&$(x^2-7^2)(x^4+30x^2+7^4)$\\ \hline
   13&$x^2-10x+13^2$& $(x^2-10x+13^2)(x^4+62x^2+13^4)$\\ \hline
   17&$x^2+30x+17^2$& $(x^2+30x+17^2)(x^2-10x+17^2)^2$\\ \hline
\end{tabular}

(Table 1)
\end{center}
\subsection{Determinants}

\begin{lemma}
{Let $\sigma_1, \sigma_2 : G_{\Q} \rightarrow \Z_2^\times$ be two
characters unramified outside $2$. If they agree on $\Frob_p$ for
  $p=3$ and $p=13$, then they are equal.}
\end{lemma}

\begin{proof}
  We follow the same argument as in the proof of lemma 5.2 in
  \cite{lly05}. Notice that the only quadratic extensions of $\Q$
  which are unramified outside 2 are $\Q(i)$, $\Q(\sqrt{2})$, and
  $\Q(\sqrt{-2})$, and
  $p=3$ is inert in the first two fields, while $p = 13$ is inert in the third.
  Hence the
  assumption implies that $\sigma_1 = \sigma_2$.
\end{proof}

Now apply the above lemma to the determinants of
$\rho_{2,l}(\Frob_p)$ and $\rho_{4,l}(\Frob_p)$ for $p=3,13$ which
can be read off from Table 1. {Since the determinant is independent
of the choice of $l$, we conclude}

\begin{cor}\label{cor:det=-p^6}
  For all odd primes $p \ne l$, we have
  $$\det(\rho_{2,l}(\Frob_p))=\chi_{-4}(p)p^2 \quad {\rm and}
  \quad
  \det(\rho_{4,l}(\Frob_p))=\chi_{-4}(p)p^6,$$
  where $\chi_{-4}(p)=\left ( \frac{ -4}{p} \right )$ is the character attached
  to the extension $\Q(i)$.
\end{cor}

\section{The case $n = 2$}
The space $S_3(\G_2)$ is 1-dimensional, spanned by
$h_2=\sqrt{E_1E_2}$.  It was shown in \cite{b-s} (corresponding to
the $\mathcal A(2)$ surface) and \cite{lly05} that Scholl's
 $l$-adic representation $\rho_{2,l}$ is modular, arising
 from a weight-3 newform $g_2=\eta(4z)^6$ with complex multiplication,
 resulting from the fact that $\E_2$ is a K3 surface with Picard number 20.
  Here $\eta(z)$ is the Dedekind eta function. It then follows from
 Scholl's theorem that $h_2$ satisfies a 3-term
 ASD congruence with $A(p)$ and $B(p)$ coming from $g_2$ since $\rho_{2,l}$ has
 degree 2.

\section{Decomposition of spaces for the case $n=4$}
\subsection{Weight-3 modular forms of $\G_4$}
Note that $\G_4$ is an index-2 subgroup of $\G_2$ with  $t_4 =
\sqrt[4]{\frac{E_1}{E_2}}$. We know $S_3(\G_4)=\{h_1,h_2,h_3\}$,
where
$$ h_1=\frac{E_1}{t_4}, \quad
h_2=\sqrt{E_1E_2}, \quad h_3=E_2\cdot t_4.$$

\subsection{Actions}

The operators $A$ and $\zeta$ defined in section 2.2 induce
actions on the cohomology space $W_{4,l}$, denoted by $A^*$ and
$\zeta^*$ respectively, and on the space of cusp forms
$S_3(\Gamma_4)$, again denoted by $A$ and $\zeta$. We examine
these actions.

 As discussed before,
\begin{eqnarray*}
 A(t_4)= t_4\big|_A&=& \frac{\omega_8}{t_4}.
\end{eqnarray*}
Therefore, the action of $A$ is defined over
$\Q(\omega_8)=\Q(i,\sqrt{2}).$ It acts on $h_1,h_2,h_3$ as  follows:
\begin{eqnarray}\label{eq:Actiononh_j}
 A\cdot h_j= h_j\big|_{A}&=& {\omega_8^{4-j}}\cdot h_{4-j}, \quad
  {\rm for} ~~j=1,2,3.
\end{eqnarray}

On the surface $\E_4$, $A$ induces the map given by
 $$A:(x,y,t_4)\mapsto (-x,iy,\frac{\omega_8}{t_4}).$$
It follows that $A^2:(x,y,t_4) \mapsto (x,-y,t_4),$  which maps
every point $P=(x,y)$ on the fibre at $t_4$ to its additive
inverse $-P=(x,-y)$. Since $A^2$ is the identity on $X_n$ and is
$-1$ on the sheaf $\CF_l$,  the action of $(A^*)^2$ on $W_{4,l}$
is multiplication by $-1$.

The map $\zeta : t_4 \mapsto \omega_8^{-2}\cdot t_4$ induces an
action on $h_j$ as follows:
\begin{equation}\label{eq:zeta_actiononh_j}
\zeta \cdot h_j=\omega _8^{2j} \cdot h_j \quad \text{for} ~j=1,2,3.
\end{equation} \noindent On the surface $\E_4$, it acts via
$$\zeta : (x,y,t_4) \mapsto (x,y, \omega_8^{-2}\cdot t_4).$$

\subsection{Decomposition of spaces}
 With $4$ and $l$ fixed, write $W$ for $W_{4,l}$.  Following Scholl
\cite{sch85b}, consider the associated $p$-adic Scholl space $V$
which contains $S_3(\G_4, \Q_p)$ as a subspace and $S_3(\G_4,
\Q_p)^{\vee}$ as a quotient.  The map $\zeta ^2$ sends $(x,y,t_4)$
to $(x,y,-t_4),$ which is defined over $\Q$ and of order 2 on $V$,
$W$. Denote by $V_{\pm}$ the eigenspaces of $\zeta ^2$ on $V$ with
eigenvalues $\pm1$. It is easy to verify that $V_-$, which
contains $h_1$ and $h_3$, is a 4-dimensional vector space over
$\Q_p$. The space $V_+$, which contains $h_2$, is fixed by
$\zeta^2$. It is the $p$-adic space attached to $S_3(\G_2)$. The
map $A$ acts on $V$ with $A^2$ being multiplication by $-1$.

On the $l$-adic side, $W$ decomposes similarly into a
2-dimensional $\Q_l$ space $W_+$ and a 4-dimensional space  $W_-$.
Likewise, $(\zeta^*)^2$ acts on $W_\pm$ as multiplication by $\pm
1$. We denote by $\rho_{\pm,l}$, or $\rho_{\pm}$ if there is no
ambiguity, the representation $\rho_{4,l}$ of $G_\Q$ restricted to
$W_{\pm}$.
\smallskip

Corollary \ref{cor:det=-p^6} implies

\begin{cor}
For all odd $p \ne l$, we have
\begin{eqnarray}\label{detrho-}
\det \rho_+(\Frob_p) = (\frac{-4}{p})p^2 \quad {\rm and} \quad
 \det \rho_-(\Frob_p) = p^4.
\end{eqnarray}
\end{cor}

\begin{lemma}\label{zeta&A}
  $\zeta^*  A^* \zeta^*  =A^*$ on $W$ and $\zeta  A \zeta  =A$ on
  $V$.
\end{lemma}
\begin{proof} On $W$ this follows from
their actions on the surface level:
\begin{multline*}
\zeta  A \zeta  (x,y,t_4)=\zeta  A (x,y,-it_4)=\zeta  (-x,iy,
\frac{\omega_8}{-it_4})\\ =(-x,iy,\frac{\omega_8}{t_4})=A(x,y,t_4).
\end{multline*}
On $V$ this follows from the actions of $\zeta$ and $A$ on
$S_3(\G_4)$.
\end{proof}

Hence $\zeta^*$ and $A^*$ (resp. $\zeta$ and $A$) generate a copy
of the quaternion group $H_8$ acting on the space $W_-$ (resp.
$V_-$).\smallskip

As we shall be comparing the actions of the Frobenius element at
 $p$ ~on $V$ and $W$, we write $F$ for its action on
 $V$ and $F_p$ for its action on $W$,
  keeping in mind that on the $p$-adic
 side there is only one Frobenius action, while on the $l$-adic
 side there are plenty. In general, for an operator $B$ acting on
 a finite-dimensional vector space $X$, denote by $\text{Char}(X, B)(T)$
 the characteristic polynomial of $B$ in variable $T$. As shown in \cite{sch85b},
 the congruences at $p$ will follow from equality of the
 characteristic
 polynomials of $F$ and $F_p$ for the relevant subspaces. Our argument
 is a refinement and generalization of that used in \cite{lly05}.

 To begin with, Scholl's theorem in \cite{sch85b} gives
 $$\text{Char}(W, F_p)(T) = \text{Char}(V, F)(T).$$
Applying the argument following \cite[4.4]{sch85b} to the matrix
$\zeta^2$, we obtain
\begin{equation}
  \text{Char}( W_+, F_p)(T)=\text{Char}(V_+, F)(T)\in \Z[T],
\end{equation} %and therefore
\begin{equation}\label{equalchar}
  \text{Char}( W_-, F_p)(T)=\text{Char}(V_-, F)(T)\in \Z[T].
\end{equation}

Since the congruences resulting from the 1-eigenspaces are for
$h_2\in S_3(\G_2)$,  whose congruence relations have been
established in \cite{lly05} and reviewed in the previous section, we
shall concentrate on the $(-1)$-eigenspaces $V_-$ and $W_-$.

Under the action of $\zeta$, the space $V_-$ further decomposes into
eigenspaces $V_{-,\pm i}$ with eigenvalues $\pm i$ respectively. In
particular, $h_1\in V_{-,i}$ and $h_3\in V_{-,-i}$.

\begin{cor}
  The constant terms of $\text{Char}(V_-,F)(T)$ and $\text{Char}(W_-,F_p)(T)$  are
  $p^4$.
\end{cor}

\section{Factorizing local $L$-factors}
In this section we shall confine ourselves to the spaces $V_-$ and
$W_-$.  Recall that $\zeta^2$, $A^2$ on $V_-$ and $(\zeta^*)^2$,
$(A^*)^2$ on $W_-$ all act as multiplication by $-1$.

\subsection{Operators on $V_-$ and $W_-$}
 In addition to $\zeta$ and $A$,
consider also the operators $$ {B_{-2}} :=(1+\zeta)A \quad \text
{and} \quad {B_{2}}: = (1-\zeta)A$$ on $V_-$. Likewise, we
introduce the operators
$$ B_{-2}^* :=A^*(1+\zeta^*) \quad \text {and} \quad B_{2}^*:
= A^*(1-\zeta^*)$$ on $W_-$. It is straightforward to check that
$$({B_{-2}})^2 = -2I = ({B_{2}})^2 \quad \text {and} \quad
(B_{-2}^*)^2 = -2I = (B_{2}^*)^2.$$ Lemma \ref{zeta&A} implies
$${B_{-2}}{B_{2}} =-{B_{2}}{B_{-2}} \quad \text {and} \quad
B_{-2}^*B_{2}^* =-B_{2}^*B_{-2}^*.$$
\begin{lemma} \label{property-pmod8} Let $p$ be an odd
prime not equal to $l$. On $W_-$ we have
\begin{itemize}
\item[(1)]
\begin{eqnarray*}
\ \zeta^*F_p=F_p (\zeta^*)^p = \left(\frac{-1}{p}\right )F_p
\zeta^*, &&
  \A^* F_p= F_p A^* (\zeta^*)^{(p-1)/2},
  \end{eqnarray*}
\item[(2)]  $$ \ B_{\pm 2}^*F_p=\left (\frac{\pm 2}{p}\right )
F_pB_{\pm 2}^*.$$
\end{itemize}
\end{lemma}
\begin{proof}\begin{itemize}
 \item[(1)] Since $F_p$ on $W_-$ is the geometric Frobenius in the Galois group, the action
can be computed as the pullback via Frobenius in the reduction
$\mod p$. Hence we examine the actions modulo $p$ on the elliptic
surface $\E_4$:

\begin{eqnarray*}
\Frob_p \zeta(x,y,t_4)&=&\Frob_p(x,y,-\omega_4
t_4)=(x^p,y^p,(-\omega_4)^p t_4^p)\\
&=&\zeta^p\Frob_p(x,y,t_4),
\end{eqnarray*}
\begin{eqnarray*}
\Frob_p A(x,y,t_4)&=&\Frob_p(-x,iy,\frac{\omega_8}{
t_4})=(-x^p,i^py^p,\frac{\omega_8^p}{ t_4^p}),
\end{eqnarray*}
 \begin{eqnarray*}
 \zeta ^{(p-1)/2}A\Frob_p(x,y,t_4)&=&\zeta ^{(p-1)/2}A(x^p,y^p, t_4^p)
 =\zeta ^{(p-1)/2}(-x^p,iy^p,\frac{\omega_8}{
t_4^p})
\\&=&\zeta ^{(p-1)/2}(-x^p,iy^p,\omega_8^{1-p} \frac{\omega_8^p}{ t_4^p})
%\\&=&\zeta ^{(p-1)/2}(-x^p,iy^p,\omega_4^{(1-p)/2} \frac{\omega_8^p}{ t_4^p})\\
\\&=&(-x^p,iy^p,\omega_4^{1-p} \frac{\omega_8^p}{ t_4^p}).
\end{eqnarray*} When $p \equiv 1 \mod
4$, this gives $\Frob_p A=\zeta ^{(p-1)/2}A\Frob_p$. When $p
\equiv 3 \mod 4$, we have $\Frob_p A=A^2\zeta^2(\zeta
^{(p-1)/2}A\Frob_p)$. On the cohomology space, the order of the
operators are reversed. Since the actions of $(A^*)^2$ and
$(\zeta^*)^2$ on
 $W_-$ are both $-1$, we obtain the desired conclusion.
\medskip
\item[(2)] By a straightforward computation using (1), we have
\begin{eqnarray*}
B_{-2}^*F_p
&=&A^*(1+\zeta^*)F_p=A^*F_p(1+(\zeta^*)^p)\\&=&F_pA^*(\zeta^*)^{(p-1)/2}(1+(\zeta^*)^p)\\
&=&\left (\frac{- 2}{p}\right ) F_pA^*(1+\zeta^*)= \left (\frac{-
2}{p}\right ) F_pB_{-2}^*. \end{eqnarray*}
 The other statement is
proved in a similar way.
\end{itemize}
\end{proof}

On the $p$-adic side, we have $\zeta h_j= \omega_4^j h_j$ and
$A(h_j)= \omega_8^{4-j}h_{4-j}$. Thus $\zeta$ on $V_-$ is defined
over $\Q_p$ whenever $\sqrt {-1} \in \Q_p$, that is, $p \equiv 1
\mod 4$, and $A$ on $V_-$ is defined over $\Q_p$ for $p \equiv 1
\mod 8$. Since $B_{2}(h_j)=(1-\zeta)\omega_8^{4-j} h_{4-j}=
(1-\omega_4^{4-j}) \omega_8^{4-j} h_{4-j}=(j-2) \sqrt{2} h_{4-j},$
hence $B_2$ on $V_-$ is defined over $\Q_p$ whenever $\sqrt{2}\in
\Q_p$, in particular, when $p \equiv 7 \mod 8$. Similarly,
$B_{-2}$ on $V_-$ is defined over $\Q_p$ whenever $\sqrt {-2} \in
\Q_p$, in particular, $p \equiv 3 \mod 8$. We record this
discussion in

\begin{lemma}\label{BF=FB}
On $V_-$ we have $\zeta F = F\zeta$ when $\sqrt {-1} \in \Q_p$,
$AF = FA$ when $\omega_8 \in \Q_p$, $B_{2}F = FB_{2}$ when $\sqrt
2 \in \Q_p$, and $B_{-2}F = FB_{-2}$ when $\sqrt {-2} \in \Q_p$.
\end{lemma}

\subsection{Factorizing local $L$-factors}

The aim of this subsection is to factor, for each odd prime $p$, the
characteristic polynomial $\text{Char}(W_-, F_p)(T)$ as a product of
two quadratic characteristic polynomials arising from a suitable
restriction of $\rho_-$.

\begin{prop}\label{uniquechar}
Let $\delta \in \{-1, -2, 2 \}$ and let $\sigma_1$ and $\sigma_2$ be
two 1-dimensional representations of $G_{\Q(\sqrt \delta)}$ over a
totally ramified extension $F$ of $\Q_2$, unramified outside the
place dividing 2.  Then $\sigma_1 = \sigma_2$ if they agree at
$\Frob_{\pi}$ for $\pi$ dividing $ 3, 13$ if $\delta = -1$ or $-2$,
and for $\pi$ dividing $ 3,  7,  13$ if $\delta = 2$.
\end{prop}

\begin{proof} Consider $\sigma=\sigma_1 /\sigma_2: G_{\Q(\sqrt \delta)}
\rightarrow F^{\times}$. Since $F^{\times}$ is a pro-2-group, so
if $\sigma\neq 1$ then its fixed field contains a quadratic
extension $K$ over $\Q(\sqrt{\d})$ unramified away from 2.

For $\delta = -1$, the possible fields $K$ are $\Q(i, \sqrt {1+i})$,
$\Q(i, \sqrt {1-i})$ in which $\pi = 3$  is inert, and $\Q(i, \sqrt
2)$ in which $\pi =  3 + 2i$ (dividing 13)  is inert.

For $\delta = -2$, the possible fields $K$ are $\Q(\sqrt {-2}, i)$,
$\Q(i\sqrt[4]{-2})$ in which $\pi = 1 +\sqrt{-2}$ (dividing 3) is
inert, and $\Q(\sqrt[4]{-2})$ in which $\pi = $   13 is inert.

For $\delta = 2$, let $\varepsilon= \sqrt
{\sqrt 2 - 1}$.  The possibilities   are  $$\Q(\sqrt 2,\varepsilon\sqrt[4] 2), \quad \Q(\sqrt 2, i \varepsilon\sqrt[4] 2), \quad \Q(\sqrt 2, \varepsilon), \quad \Q(\sqrt 2, i
\varepsilon)$$ in which $\pi = 3$ is inert, $\Q(\sqrt[4] 2)$,
$\Q(i \sqrt[4] 2)$ in which $\pi = 13$ is inert, and $\Q(\sqrt {-2},
i)$ in which $\pi = 3 + \sqrt 2$ (dividing 7)
 is inert.

\end{proof}

We proceed to choose a quadratic character of $G_{\Q(\sqrt
\delta)}$ unramified outside the unique place $v$ dividing 2 which
will be needed for our purpose. Denote by $\theta_{\delta}$ the
quadratic character attached to the extension $\Q(i, \sqrt 2)$
over $\Q(\sqrt \delta)$. When $\delta = -1$, the primes $p \equiv
3, 7 \mod 8$ are inert in $\Q(i)$ with residue field $\mathbb
F_{p^2}$ containing a square root of 2, hence they split in $\Q(i,
\sqrt 2)$. The primes $p \equiv 1 \mod 8$ split in $\Q(i)$ with
residue field $\mathbb F_p$ containing a square root of $2$, hence
these places also split in $\Q(i, \sqrt 2)$. We have
$\theta_{-1}(\Frob_v) = 1$ at places $v$ dividing $p \equiv 1 \mod
8$. The primes $p \equiv 5 \mod 8$ split in $\Q(i)$ with residue
field $\mathbb F_p$ in which $2$ is not a square, thus
$\theta_{-1}(\Frob_v) = -1$ only at places $v$ above $p \equiv 5
\mod 8$. Similarly, $\theta_{-2}(\Frob_v) = -1$ only at places $v$
dividing $p \equiv 3 \mod 8$ and $\theta_2(\Frob_v) = -1$ only at
places $v$ dividing $p \equiv 7 \mod 8$.
\smallskip

Write $B_{-1}^*$ for $\zeta$. Recall that $B_{-1}^*$, $B_{-2}^*$
and $B_{2}^*$ act on the 4-dimensional space $W_-$ over $\Q_l$,
and satisfy $(B_{-1}^*)^2 = -I$, $(B_{-2}^*)^2 =(B_{2}^*)^2 =
-2I$.
 The commuting relations between these
operators and the Frobenius elements described in Lemma
\ref{property-pmod8} show that for $\delta = -1, 2, -2$,
$B_\delta^*$ is defined over $\Q(\sqrt{\delta})$. For each
$\delta$, choose a prime $l$ so that $\Q_l(B_\delta^*)$ is a
quadratic extension of $\Q_l$. Regarding
 $W_-$ as a 2-dimensional vector space over
 $\Q_l(B_\delta^*)$,
  we may restrict the representation $\rho_-$ to
 $H_\d := G_{\Q(\sqrt{\delta})}$, obtaining a 2-dimensional representation $\rho_{-,\delta}$ of
 $H_\d$ on the space $W_-$ over $\Q_l(B_\delta^*)$.

\begin{theorem}\label{Lofrho-}
The $L$-factor of the 4-dimensional representation $\rho_-$ at an
odd prime $p \ne l$ is equal to the product of the $L$-factors of
the 2-dimensional representation $\rho_{-,\delta}$ over the places
of $\Q(\sqrt \delta)$ dividing $p$ for $\delta = -1, -2, 2$.
\end{theorem}
\begin{proof}
Extend scalars to $\overline{\Q}_l$. Let $J_{-1}=B_{-1}^*$,
$J_{-2}=B_{-2}^*/\sqrt{2},$ and $J_{2}=B_{2}^*/\sqrt{2}$ as
endomorphisms of $W_-$. By Lemmas \ref{zeta&A} and
\ref{property-pmod8}, we have
\begin{itemize}
  \item[(i)] $ J_\delta^2 = -I$ {\rm ~for} $~\delta = -1, 2, -2,$ and $J_{-1}J_2=J_{-2}=-J_2J_{-1}$.
  \item[(ii)] If $\varepsilon_\d: G_{\Q} \rightarrow \{\pm 1\}$ are
  the characters corresponding to the fields $\Q(\sqrt{\delta})$,
  then $\rho_-(g)J_{\d}=\varepsilon_{\d}(g)J_{\d} \rho_-(g)$ for all $g\in
  G_{\Q}$ and $\d=-1,2,-2$.
\end{itemize}
So the $J_\d$ generate a subalgebra of $\textrm{End}(W_-)$
isomorphic to $M_2(\overline{\Q}_l)$. Hence  one can find a basis
of $W_-$ with respect to which the $J_\d$ are represented by
\begin{equation}
  J_{-1}=\begin{pmatrix}
    iI_2&0\\0&-iI_2
  \end{pmatrix}, \quad J_2=\begin{pmatrix}
    0&I_2\\-I_2&0
  \end{pmatrix}, \quad J_{-2}=\begin{pmatrix}
   0& iI_2\\iI_2&0
  \end{pmatrix}.
\end{equation}
For $g \in H_\d$, $\rho_-(g)$ commutes with $J_\d$, hence is of
the form
\begin{equation*}
  \begin{pmatrix}
    P&0\\0&S
  \end{pmatrix}, \quad \begin{pmatrix}
    P&Q\\-Q&P
  \end{pmatrix}, \quad \begin{pmatrix}
    P&Q\\Q&P
  \end{pmatrix},
\end{equation*} according as $\d =-1,2,-2$. Setting
$N=\bigcap_{\d}
H_\d=G_{\Q(i,\sqrt{2})}$, we get
\begin{eqnarray*}
  \rho_-(N)=\left \{ \begin{pmatrix}
    P&0\\0&P
  \end{pmatrix} \right \},&& \rho_-(H_{-1}\setminus N)=\left \{\begin{pmatrix}
    P&0\\0&-P
  \end{pmatrix} \right \},\\
  \rho_-(H_2\setminus N)=\left \{\begin{pmatrix}
    0&Q\\-Q&0
  \end{pmatrix} \right \},&& \rho_-(H_{-2}\setminus N)=\left \{\begin{pmatrix}
    0&Q\\Q&0
  \end{pmatrix} \right \}.
\end{eqnarray*}
Let $\sigma_{-1}: H_{-1} \rightarrow GL_2(\overline{\Q}_l)$ be the
representation mapping $g \in H_{-1}$ to the matrix $P$ in the
expression of $\rho_-(g)$ above. Identifying the character
$\theta_{-1}$ of $\Q(i,\sqrt{2})/\Q(i)$ discussed above with the
character on $H_{-1}/N$, we have
\begin{equation*}
  \rho_-|_{H_{-1}}=\begin{pmatrix}{
    \sigma_{-1}}&0\\0&{
    \sigma_{-1}} \otimes \theta_{-1}
  \end{pmatrix}=\text{Ind}_{H_{-1}}^{G_{\Q}}(\sigma_{-1})|_{H_{-1}}.
\end{equation*}
This shows that $\sigma_{-1}$ is $\rho_{-,-1}$.

For $\d = \pm 2$ we can choose a basis so that $J_\d$ is
represented by $\begin{pmatrix}
    iI_2&0\\0&-iI_2
  \end{pmatrix}$, $J_{-1}$ by $\begin{pmatrix}
    0&I_2\\-I_2&0
  \end{pmatrix}$ and the third matrix determined by property (i).
  A similar argument then shows that
  \begin{equation}\label{induced}
  \rho_-|_{H_{\d}} =\begin{pmatrix}{
    \rho_{-, \d}}&0\\0&{
    \rho_{-, \d}} \otimes \theta_{\d}
  \end{pmatrix} = \text{Ind}_{H_{\d}}^{G_{\Q}}(\rho_{-,
  \d})|_{H_{\d}} \quad \text{ for} ~~\d = 2, -2 ~\text{(and} ~-1).
  \end{equation}
Therefore the local $L$-factors attached to $\rho_-$ have the
asserted property.
\end{proof}
\smallskip

\subsection{The determinants of the restricted representations}
 Using Table 1 in Sec.
\ref{subsec:trace}, we obtain the following:
\begin{eqnarray}\label{charpolyat357}
\begin{split}
\text{Char}(W_-, F_3)(T)&=(T^2 - 2 \sqrt {-2}T - 3^2 )(T^2 +
2\sqrt{-2}T - 3^2),\\
\text{Char}(W_-, F_7)(T)&= (T^2 - 8 \sqrt{-2}T - 7^2)(T^2 + 8
\sqrt{-2}T -
7^2),\\
\text{Char}(W_-, F_{13})(T)&= (T^2+20iT-13^2)(T^2-20iT-13^2).\\
\end{split}
\end{eqnarray}
As this is the unique way to factor $\text{Char}(W_-, F_p)(T)$ for
$p = 3, 7, 13$
 into a product of two degree two
polynomials with opposite coefficients for $T$, applying Theorem
\ref{Lofrho-}, we obtain the following table of the values of
$\det \rho_{-,\delta}(\Frob_v)$:

\begin{center}
\begin{tabular}{|c|c|c|c|}
\hline $\det \rho_{-,\delta}(\Frob_v)$&$v$ \text{ above } 3& $v$
\text{ above } 7& $v$
\text{ above { 13}}\\
\hline
$\delta=-1$&$3^4$&$7^4$&$-13^2$\\
$\delta=-2$&$-3^2$&$7^4$&$13^4$\\
$\delta=2$&$3^4$&$-7^2$&$13^4$\\
\hline \end{tabular}
\end{center}

\noindent Combined with Proposition \ref{uniquechar}, we conclude

\begin{cor}\label{detrhodelta}
 For $\delta = -1, -2, 2$, at a place $v$ of $\Q(\sqrt \delta)$
 dividing
 an odd prime $p \ne l$, we have $$\det \rho_{-,\delta}(\Frob_v)
 = \theta_{\delta}(\Frob_v) Nv^2,$$ where $Nv$ is the norm
 of $v$.
\end{cor}

As remarked earlier, for each odd prime $p \ne l$ and not $\equiv 1
\mod 8$, there is a quadratic extension $\Q(\sqrt \delta)$ in which
$p$ splits and $\theta_{\delta}(\Frob_v) = -1$ at places $v$
dividing $p$, hence we can combine Theorem \ref{Lofrho-} with the
above corollary to give a more detailed description of the
factorization of local factors of $\rho_-$.

\begin{cor}\label{roughdet}
For each odd prime $p \ne l$, there is a constant $a_p$, not
depending on $l$, such that $\text{Char}(W_-, F_p)(T) = (T^2 -
a_pT + p^2)^2$ if $p \equiv 1 \mod 8$, and $\text{Char}(W_-,
F_p)(T)=(T^2 - a_pT -p^2)(T^2 + a_pT -p^2)$ otherwise.
\end{cor}

\section{ASD congruences for $S_3(\G_4)$}
 A cusp form $h(z)$ in $S_3(\G_4)$ with $2$-integral Fourier coefficients
$a(n)$ is said to satisfy a 3-term Atkin-Swinnerton-Dyer
congruence relation at an odd prime $p$ if there exist algebraic
integers $A_p$ and $B_p$ with $|\sigma(A_p)| \le 2 p$ and
$|\sigma(B_p)| = p^2$ for all embeddings $\sigma$ so that for all
$n \ge 1,$
\begin{eqnarray}\label{e:0.5}
\frac{a(np) - A_pa(n) + B_p a(n/p)}{(pn)^2}
\end{eqnarray}
is integral at some place above $p$.
 For brevity, we refer to this as $h$ satisfying an ASD congruence
 at $p$ given by $T^2-A_pT+B_p$.

 In this section we shall prove the following two main results.

\begin{theorem}\label{charpolyoffrob}
For an odd prime $p \ne l$, $\text{Char}(W_-, F_p)(T)$ has the
following factorization for some $A_p \in \Z$.
\begin{enumerate}
\item If $p \equiv 1 \mod 8$, $\text{Char}(W_-, F_p)(T) = (T^2 -
A_pT + p^2)^2$; \item If $p \equiv 5 \mod 8$, $\text{Char}(W_-,
F_p)(T) = (T^2 - iA_pT - p^2)(T^2 + iA_pT - p^2)$; \item If $p
\equiv 3 \ {\rm or } \ 7 \mod 8$,
$$\text{Char}(W_-, F_p)(T) = (T^2 - \sqrt{-2}A_pT - p^2)(T^2 +
\sqrt{-2}A_pT - p^2).$$
\end{enumerate}
\end{theorem}
\smallskip

 \begin{theorem}[ASD congruence for $S_3(\G_4)$]\label{ASDcongruenceconclusion}
The ASD congruence holds on the space $S_3(\G_4) = <h_1, h_2, h_3>$.
More precisely, $h_2$ lies in $S_3(\G_2)$ and it satisfies the ASD
congruence relations with the congruence form $g_2(z)= \eta(4z)^6$.
For each odd prime $p$, the subspace $<h_1, h_3>$ has a basis
depending on the residue of $p$ modulo $8$ satisfying a 3-term ASD
congruence at $p$ as follows.
\begin{enumerate}
\item If $p \equiv 1 \mod 8$, then both $h_1$ and $h_3$ satisfy
the 3-term ASD congruence at $p$ given by $T^2 - A_pT + p^2$;
\item If $p \equiv 5 \mod 8$, then $h_1$ (resp. $h_3$)  satisfies
the 3-term ASD-congruence at $p$
 given by $T^2 - iA_pT - p^2$ (resp. $T^2 + iA_pT - p^2$);
 \item If $p \equiv 3 \mod 8$ (resp. $p \equiv 7 \mod 8$), then $h_1\pm  h_3$
 (resp. $h_1\pm i h_3$) satisfy the 3-term
 ASD congruence at $p$ given by
  $T^2 \mp \sqrt{-2}A_pT -p^2$, respectively.
  \end{enumerate}
  The polynomials above are factors of $\text{Char}(W_-,
  F_p)(T)$ as shown in Theorem \ref{charpolyoffrob}.
  \end{theorem}
\smallskip

\begin{cor}\label{charpolymod2} For all primes $p$ we have
$$ \text{Char}(W_-, F_p)(T) \equiv T^4 + 1 \mod 2.$$
\end{cor}
\begin{proof}
  When $p \equiv 3 \mod 4$, the verification is straightforward. When $p \equiv 1 \mod
  8$,   $\Tr\rho_-(\Frob_p)=2A_p \equiv 0 \mod 4$ by Lemma
\ref{lem:ap=0mod8ifp=1mod8}, hence $A_p$ is even. When $p \equiv 5
  \mod 8$, Lemma \ref{lem:ap=0mod8ifp=1mod8} implies $\Tr
  \rho_-(\Frob_p^2)=2(2p^2-A_p^2) \equiv 0 \mod 4$, thus $A_p$ is
  even.
\end{proof}

\subsection{Proof of Theorems  \ref{charpolyoffrob} and \ref{ASDcongruenceconclusion}}
The space $S_3(\G_4)$ is spanned by $h_1, h_2$, and $h_3$. We know
that $h_2$ generates the space $S_3(\G_2)$ and it satisfies ASD
congruence relations as proved in \cite{lly05} and reviewed in
section 4. So it remains to prove the theorem for the subspace
$<h_1, h_3>$ as stated. The general strategy is to find suitable
operators $B$ and $B^*$ acting on $V_-$ and $W_-$ respectively,
commuting with the action of the Frobenius at $p$, such that the
characteristic polynomials of $F$ and $F_p$ on the respective
eigenspace of $B$ and $B^*$ with the same eigenvalue agree. Since
the characteristic polynomial of $F_p$ is independent of the
choice of $l \ne p$, we shall choose $l \equiv p \mod 8$ so that
$\Q_l$ always contains the eigenvalues of $B^*$.

\subsubsection{Case I. $p \equiv 1
\mod 8$} The eigenspaces $V_{-,\pm i}$ of $\zeta$ on $V_-$ (resp.
 $W_{-,\pm i}$ of $\zeta^*$ on $W_-$) with
eigenvalue $\pm i$ are $F$- (resp. $F_p$-) invariant. Further,
since $A^*$ commutes with the action of $F_p$ and it maps
$W_{-,i}$ to $W_{-,-i}$ isomorphically, we get
$$\text{Char}(W_{-,i},F_p)(T)=\text{Char}(W_{-,-i},F_p)(T)=T^2-A_pT+B_p$$
for some constants $A_p$ and $B_p$. That $A_p \in \Z$ and $B_p =
p^2$ follows from (\ref{induced}) and Corollary \ref{roughdet}.

By \cite{sch85b}, we know that $F_p$ on $W_-$ and $F$ on $V_-$
have the same characteristic polynomials, and the same is true for
$\zeta^* F_p$ and $\zeta F$. This implies
$$\text{Char}(W_{-,\pm i},F_p)(T)=\text{Char}(V_{-,\pm i},F)(T).$$
Combined with $h_1 \in V_{-,i}$ and $h_3 \in V_{-,-i}$, this
proves the asserted ASD-congruence.
\smallskip

To prove the remaining cases, we shall need the following Lemma.
Let $\d=-1, 2, -2$ and $B_{\d}^*$ be as in the previous section so
that $(B_{\d}^*)^2 = \l$, where $\l = -1, -2, -2$ according as
 $\d = -1, 2, -2$. Let  $p$  and $l$ be primes such that $\Q_p$ and  $\Q_l$
contain $\sqrt{\d}$. Then $W_-$ (resp. $V_-$) decomposes into a
direct sum of eigenspaces $W_{-,\pm \sqrt{\l}}$ (resp. $V_{-,\pm
\sqrt{\l}}$) of $B_{\d}^*$ (resp. $B_{\d}$), which are invariant
under the action of the Frobenius at $p$ by Lemmas
\ref{property-pmod8} and \ref{BF=FB}.

 \begin{lemma}\label{lem:7.4} With the above notation, if
\begin{equation}\label{equalcharpoly}
  \text{Char}(W_{-,+ \sqrt{\l}},
  B_{\d}^* F_p)(T)=\text{Char}(W_{-,- \sqrt{\l}},
  B_{\d}^* F_p)(T)=T^2-a_pT+b_p
  \end{equation} for some constants $a_p$ and $b_p$,  then
  $a_p$ and $b_p$ lie in $\l\Z$, and  $$\text{Char}(W_{-,\pm \sqrt{\l}},
  F_p)(T)=\text{Char}(V_{-,\pm \sqrt{\l}},
  F)(T).$$
\end{lemma}
\begin{proof} Write $B_{\d}=\sum k_i A_i$ as a linear combination of $A_i\in \SL$
with coefficients $k_i\in \Z$. The traces of $B_{\d}F$ on $V$ and
$V_+$ are equal to $\sum  k_i \Tr (A_iF)$ on the respective
spaces, and likewise for $B_{\d}^*F_p$ on $W$ and $W_+$. We
conclude from Sec. 4.4 of \cite{sch85b} that the traces of
$B_{\d}F$ on $V_-$ and $B_{\d}^*F_p$ on $W_-$ agree and they are
in $\Z$.

  On the $p$-adic side, using $\zeta(h_j)= \omega_4^j h_j$ and
$A(h_j)= \omega_8^{4-j}h_{4-j}$,  one finds that each $(\pm
\sqrt{\l})$-eigenspace of $B_{\d}$ on $V_-$ contains a linear
combination $h_{\pm \l}$ of $h_1$ and $h_3$; further, the dual of
$h_{\mp \l}$ appears in the quotient of $(\pm
\sqrt{\l})$-eigenspace modulo $h_{\pm \l}$. This shows that each
eigenspace of $B_{\d}$ on $V_-$ is 2-dimensional.

   Let $\sqrt{\l}\al_1$ and $\sqrt{\l} \al_2$ (resp. $-\sqrt{\l}
  \al_3$ and $-\sqrt{\l} \al_4$) be the eigenvalues of $B_{\d}^*F_p$ on
  $W_{-,\sqrt {\l}}$ (resp. $W_{-,-\sqrt {\l}}$) so that $\al_1, \al_2, \al_3,\al_4$
  are the eigenvalues of $F_p$ on $W_-$. In view of
  (\ref{equalcharpoly}), we may assume $\al_3=-\al_1$ and
  $\al_4=-\al_2$. Thus $T^2-a_pT+b_p = (T - \sqrt {\l}\al_1)(T - \sqrt
  {\l}\al_2)$ and
  $$\text{Char}(W_-, F_p)(T) = (T^2 - \al_1^2)(T^2 - \al_2^2)
 \in \Z[T].$$
Therefore $b_p= \l \al_1 \al_2 = \pm \l p^2 \in \l\Z$ and $a_p^2 =
\l(\al_1 + \al_2)^2 = \l(\al_1^2 +
 \al_2^2) + 2b_p \in \l\Z$. This combined with $2a_p = \Tr B_{\d}^*F_p \in \Z$
 implies $a_p \in \l\Z$ since $\l$
 is square-free.

It remains to prove the last assertion. Since $F$ on $V_-$ has the
same eigenvalues as $F_p$ on $W_-$, first consider the situation
that one of $\al_1, \al_2$ is an eigenvalue of $F$ on
$V_{-,+\sqrt{\l}}$. Due to the symmetry on
 $\al_1$ and $\al_2$, we may assume it is $\al_1$. Then
there are three possibilities for the second eigenvalue of $F$ on
$V_{-,+\sqrt{\l}}$~: (i) $\al_2$, (ii) $-\al_1,$ and (iii) $
  -\al_2.$ If it is case (i), then we are done. If it
  is case (ii), then the eigenvalues of $B_{\d}F$ are $\sqrt{\l}\al_1, -\sqrt{\l}\al_1,
  -\sqrt{\l}\al_2,$ and $ \sqrt{\l} \al_2$ so that $B_{\d}F$ has zero trace. As
$B_{\d}F$ and $B_{\d}^*F_p$ have the same trace, we conclude that
$\al_1=-\al_2$ and hence the assertion also holds. If it is case
(iii), then the eigenvalues for $B_{\d}F$ are $\sqrt{\l}\al_1,
-\sqrt{\l}\al_2,
  \sqrt{\l}\al_1, -\sqrt{\l}\al_2$. Since the traces of
$B_{\d}F$ and $B_{\d}^*F_p$ are the same, one concludes that
$\al_2=0$, which contradicts $|\al_2|=p$. So this case cannot
occur. Finally we note that case (ii) also takes care of the
situation that one of $-\al_1$ and $-\al_2$ is an eigenvalue of
$F$ on $V_{-,+\sqrt{\l}}$. This completes the proof of the Lemma.
\end{proof}

\subsubsection{Case II. $p \equiv 5 \mod 8$} In this case $\zeta$ is defined over $\Q_p$.
Since $\zeta^*  F_p=F_p \zeta^*$ by Lemma \ref{property-pmod8},
$F_p$ leaves invariant the eigenspaces $W_{-,i}$ and $W_{-,-i}$ of
$\zeta^*$ on $W_-$. Further, by Lemma \ref{property-pmod8}, $A^*$
commutes with $\zeta^* F_p$ and it maps $W_{-,i}$ isomorphically
to $W_{-,-i}$, therefore
$$\text{Char}(W_{-,i},\zeta^* F_p)(T)=\text{Char}(W_{-,-i},\zeta^*
F_p)(T)=T^2-a_pT+b_p$$ for some constants $a_p$ and $ b_p$. By
Lemma \ref{lem:7.4}, we have $a_p, b_p  \in \Z$ and
$$\text{Char}(W_{-,\pm i},\zeta F_p)(T)=
\text{Char}(V_{-,\pm i},\zeta F)(T).$$ Combining with
(\ref{induced}) and Corollary \ref{roughdet}, and noticing $h_1
\in V_{-,i}$ and $h_3 \in V_{-,-i}$ , we obtain the desired
assertions for this case.

\subsubsection{Case III. $p \equiv 3 \mod 8$} Denote by $V_{-,\pm
\sqrt{-2}}$ (resp. $W_{-,\pm \sqrt{-2}}$) the eigenspaces of
$B_{-2}$ on $V_-$ (resp. $B_{-2}^*$ on $W_-$) with eigenvalue $\pm
\sqrt{-2}$.
 By Lemma \ref{property-pmod8}, the eigenspaces of $B_{-2}^*$ are
invariant under $F_p$, $B_{2}^*$ commutes with $B_{-2}^*F_p$ and
it maps $W_{-, \sqrt{-2}}$ isomorphically to $W_{-, -\sqrt{-2}}$.
So there are $a_p$ and $b_p$ such that
$$\text{Char}(W_{-,\sqrt{-2}},B_{-2}^* F_p)(T)= \text{Char}(W_{-,-\sqrt{-2}},B_{-2}^* F_p)(T)=T^2-a_pT+b_p.$$
It follows from (\ref{induced}), Corollary \ref{roughdet} and
Lemma \ref{lem:7.4} that the characteristic polynomial of $F_p$
has the asserted form and $$\text{Char}(W_{-,\pm
\sqrt{-2}},B_{-2}^* F_p)(T) = \text{Char}(V_{-,\pm
\sqrt{-2}},B_{-2}^* F)(T).$$

 Finally it is a straightforward computation, using the actions of $A$ and $\zeta$ given by
(\ref{eq:Actiononh_j}) and (\ref{eq:zeta_actiononh_j}), to check
that $h_1 \pm h_3$ are eigenfunctions of ${B_{-2}}$ on $V_-$ with
eigenvalues $\pm \sqrt{-2}$, respectively.

\subsubsection{Case IV. $p \equiv 7 \mod 8$}
The detailed analysis parallels the previous case with the roles
of $B_{-2}^*$ and $B_{2}^*$ interchanged; the eigenvalues of
$B_{2}$ on $V_-$ (resp. $B_{2}^*$ on $W_-$) are $\pm \sqrt {-2}$
with eigenspaces $V_{-,\pm\sqrt {-2}}$ (resp. $W_{-,\pm\sqrt
{-2}}$). In this case one checks that $h_1\pm ih_3 \in
V_{-,\pm\sqrt {-2}}$.

This completes the proof of Theorem \ref{charpolyoffrob} and
Theorem \ref{ASDcongruenceconclusion}.

\section{Modularity of $\rho_{4,l}$}

\subsection{An automorphic representation} Let $K = \Q(i, 2^{1/4})$,
which is a Galois extension over $\Q$ with Galois group
$\Gal(K/\Q)$
  dihedral of order 8. It is a degree 4 extension over $\Q(i)$ with
  $\Gal(K/\Q(i))$ cyclic of order 4. The Artin reciprocity map from the
idele class group of $\Q(i)$ to $\Gal(K/\Q(i))$ followed by an
isomorphism from $\Gal(K/\Q(i))$   to the group $<i>$ generated by
$i \in \mathbb C^\times$ yields an idele class character of
$\Q(i)$ of order 4, denoted by $\chi$. It is ramified only at the
place $1+i$ (above 2) of $\Q(i)$. Its values at the places above
the odd primes $p$ of $\Q$ are as follows:
\begin{enumerate}
 \item  For $p \equiv 3 \mod 4$, it remains a prime in $\Q(i)$. We know
 that $2$ (resp. $-2$) is a square in $\Z/p\Z$ when $p \equiv 7$ (resp. 3) $\pmod 8$, and is thus a fourth power in the residue
  field of $\Q(i)$ at $p$. Consequently, $p$ splits completely in $K$ and we have
$\chi(p) = 1$.

\item For $p \equiv 1 \mod 8$, it splits into two places $v_1,
v_2$ of $\Q(i)$. Since 2 is a square modulo 8, we have
$\chi(v_1)=\chi(v_2) = \pm 1 \equiv 2^{(p-1)/4} \mod p$.

\item For $p \equiv 5 \mod 8$, it splits into two places $v_1,
v_2$ of $\Q(i)$. We have $\chi(v_1) = \chi(v_2)^{-1} = \pm i$, again
determined by $2^{(p-1)/4} \mod v_1 (\text{resp.  }v_2)$ since 2 is
not a square modulo $p$. The opposite sign comes from the fact that
$v_1$ and $v_2$ are complex conjugates, and $2^{(p-1)/4}$ is a
fourth root of unity, which is congruent to $i$ modulo one prime and
$-i$ modulo the other.
\end{enumerate}
Thus $\chi^2$ is a quadratic character of the idele class group of
$\Q(i)$ which is equal to $-1$ only at the places above $p \equiv 5
\mod 8$. In other words, $\chi^2 = \epsilon_{\Q(i)}$ in section 6.

Let
\begin{eqnarray}\label{eq:f-i}
\begin{split}
f_1(z)&=& \frac{\eta(2z)^{12}}{\eta(z)\eta(4z)^5}=q^{1/8}(1+q-10q^2+\cdots)=\sum_{n\ge 1}a_1(n)q^{n/8},\\
f_3(z)&=&\eta(z)^5\eta(4z)=q^{3/8}(1-5q+5q^2+\cdots)=\sum_{n\ge 1}a_3(n)q^{n/8},\\
f_5(z)&=& \frac{\eta(2z)^{12}}{\eta(z)^5\eta(4z)}=q^{5/8}(1+5q+8q^2+\cdots)=\sum_{n\ge 1}a_5(n)q^{n/8},\\
f_7(z)&=&\eta(z)\eta(4z)^5=q^{7/8}(1-q-q^2+\cdots)=\sum_{n\ge
1}a_7(n)q^{n/8},
\end{split}
\end{eqnarray}
\noindent and define \begin{equation} \label{eq:f'}f'=f'(z) = f_1(z)
+ 4f_5(z) + 2\sqrt{-2}(f_3(z) - 4f_7(z)) = \sum_{n \ge
1}a(n)q^{n/8}.\end{equation} It is easy to verify that $f'(8z)$ is a
classical cuspform of level 256, weight 3, and quadratic character
$\chi_{-4}$ associated to $\Q(i)$, and that it is an eigenform of
the Hecke operators at odd primes. The twists of $f'$ by the three
quadratic characters of $(\Z/8\Z)^{\times}$ also have the same
property.

Let $\rho'$ be the $\lambda$-adic representation of $G_{\Q}$
associated to $f'$. Then $L(s, f') = \prod_{p \ne 2}1/(1 -
a(p)p^{-s} + (\frac{-4}{p})p^{2-2s})$ is equal to $L(s, \rho')$ if
we remove the factor at $l$ divisible by $\lambda$. As the Fourier
coefficients of $f'$ lie in $\Z[\sqrt{-2}]$, the $\lambda$-adic
representation $\rho'$ yields an action of $G_\Q$ on a 2-dimensional
vector space over $\Q_l(\sqrt {-2})$.

Denote by $\rho$ the restriction of $\rho'$ to the index-2
subgroup $G_{\Q(i)}$ so that it is a degree two $\lambda$-adic
representation of $G_{\Q(i)}$. Corresponding to $\rho$ is the cusp
form $f$ for $GL_2(\Q(i))$, which is the lifting of $f'$ to a form
over $\Q(i)$ under the base change by Langlands (c.f.
\cite{langlands-basechange}). Since $f'$ and $f'$ twisted by
$\chi_{-4}$ both lift to the same form $f$, corresponding to the
representation $\rho \otimes \chi$ is the cusp form $f_{\chi}$,
called $f$ twisted by $\chi$, for $GL_2(\Q(i))$,
 whose L-function is
\begin{eqnarray*}
 &&L(s, f_{\chi})=\prod_{p \equiv 3 \mod 4} \frac{1}{(1- a(p)p^{-s}-
p^{2-2s})(1+ a(p)p^{-s} - p^{2-2s})} \\&&
 \times \prod_{p \equiv 1 \mod 8}\prod_{v_1, v_2|p} \frac{1}{(1 - \chi(v_1)a(p)p^{-s} + p^{2-2s})
 (1 - \chi(v_2)a(p)p^{-s} + p^{2-2s})}  \\
&& \times \prod_{p \equiv 5 \mod 8} \frac{1}{(1 - a(p)ip^{-s} -
p^{2-2s})(1 + a(p)ip^{-s} - p^{2-2s})},\end{eqnarray*} which is
$L(s, \rho \otimes \chi)$ after removing the factors at the places
dividing $l$. In the formula above, when $p \equiv 1 \mod 8$, $v_1$
and $v_2$ are two places of $\Q(i)$ dividing $p$, and $\chi(v_1) \chi(v_2)= \pm 1 \equiv 2^{(p-1)/4} \mod p$, as discussed above.
Observe that while there are two choices for $\chi$, the
$L$-function above is independent of the choice. Moreover, the
$L$-function remains the same if $f'$ is twisted by any quadratic
character of $(\Z/8\Z)^\times$. Finally, $\rho \otimes \chi$ can be
realized as a 2-dimensional representation of $G_{\Q(i)}$ over
$\Q_l(i)$.

Recall that the representation $\rho_{-,-1}$, the restriction of
$\rho_-$ to $G_{\Q(i)}$, can be viewed as a representation of
$G_{\Q(i)}$ to the 2-dimensional vector space $W_-$ over
$\Q_l(\zeta)$, and its associated $L$-function, written as an
Euler product over the odd primes, agrees with the $L$-function
attached to $\rho_-$, as shown in Theorem \ref{Lofrho-}.  Take
 $l = 2$ so that
we may identify $\Q_2(\zeta)$ with $\Q_2(i)$ such that $\rho_{-,-1}$
and $\rho \otimes \chi$ have the same local $L$-factors at the place
$3+2i$ and hence also at $3-2i$. This is possible from comparing
(\ref{charpolyat357}) and the Fourier coefficient $a(13)$ of $f'$.

 We want to show that $\rho_-$ and $\rho
\otimes \chi$ have the same local $L$-factors over the odd primes,
and are hence isomorphic. This will follow from

\begin{theorem}\label{thm:main} The two representations $\rho_{-,-1}$ and
$\rho \otimes \chi$ of $G_{\Q(i)}$ have isomorphic
semi-simplifications.
\end{theorem}

\begin{proof}
 Note that both representations are unramified outside
the place $1+i$. In view of Theorem \ref{charpolyoffrob}, Corollary
\ref{detrhodelta} and the explicit expression of the $L$-function
attached to $f_{\chi}$, the two representations have the same
determinants and the same local $L$-factors at $3, 3+2i, 3-2i,$ and
$7$. Moreover, by Corollary \ref{charpolymod2} and the definition of
$f'$, for both representations, the characteristic polynomials of
the Frobenius elements at places outside $1+i$ are all congruent to
$T^2 + 1$ modulo 2.

 To proceed, we use the following result of Livn\'e
\cite{liv87}, which is an extension of Serre's method \cite{ser_1} applied
to the case of representations with even trace.

\begin{theorem}[Livn\'e] Let $K$ be a global field, $S$ a finite set of
 places of $K$, and $E$ a finite extension of $\Q_2$. Denote the maximal ideal
 in the ring of integers of $E$ by $\CP$ and the compositum of all quadratic
 extensions of $K$ unramified outside $S$ by $K_S$. Suppose $\rho_1, \rho_2 :
 G_K \rightarrow GL_2(E)$ are continuous representations, unramified outside $S$,
 and furthermore satisfying
\begin{enumerate}

\item $\Tr \rho_1 \equiv \Tr \rho_2 \equiv 0 \mod \CP$ and
$\det \rho_1 \equiv \det \rho_2 \mod \CP$;

\item There exists a set $T$ of places of $K$, disjoint from $S$, for which
\begin{enumerate}
\item The image $T'$ of the set $\{\Frob_t\}_{t \in T}$ in (the $\Z/2\Z$-
vector space) $\Gal(K_S/K)$ has the property that any cubic homogeneous
polynomial on $\Gal(K_S/K)$ which vanishes on $T'$ vanishes on the vector
space $\Gal(K_S/K)$;

\item $\Tr ~\rho_1(\Frob_t) = \Tr ~\rho_2(\Frob_t)$ and $\det
\rho_1(\Frob_t) = \det \rho_2(\Frob_t)$   for all $t \in T$.
\end{enumerate}
\end{enumerate}
Then $\rho_1$ and $\rho_2$ have isomorphic semi-simplifications.
\end{theorem}

Apply the above theorem to the case $\rho_1 = \rho_{-,-1}$,
$\rho_2 = \rho \otimes \chi$ with $K = \Q(i)$, $E =
\Q_2(i,\sqrt{2})=\Q_2(\omega_8)$ and $S= \{1+i\}$. Then $K_S =
\Q(i, \sqrt 2, \sqrt {1+i})$ is a biquadratic extension of $K$
with the third quadratic extension being $\Q(i, \sqrt {1-i})$.
Choose the set $T$ to consist of the three places $3, 3+2i$ and
$3-2i$ of $\Q(i)$, which split in $\Q(i, \sqrt 2)$, $\Q(i, \sqrt
{1-i})$, and $\Q(i, \sqrt {1+i})$, respectively, and are inert in
the other two quadratic extensions of $\Q(i)$. Thus the Frobenius
elements at the places in $T$ are precisely the three nontrivial
elements of $\Gal(K_S/K)$. Further, $\rho_1$ and $\rho_2$ have the
same local $L$-factors at these three places, as observed before.
Therefore all the conditions are satisfied, and hence $\rho_1$ and
$\rho_2$ have isomorphic semi-simplifications.
\end{proof}

\begin{cor} Representations $\rho_-$ and $\rho \otimes \chi$ have the
same local $L$-factors over all odd primes $p$.
\end{cor}

Combined with Theorem \ref{charpolyoffrob}, we obtain an explicit
description of the three-term ASD congruence relation in Theorem
\ref{ASDcongruenceconclusion}.

\begin{theorem}\label{ASDexplicit} [ASD congruence for the space $<h_1, h_3>$]
For each odd prime $p$, the subspace $<h_1, h_3>$ has a basis
depending on the residue of $p \mod 8$ satisfying a 3-term ASD
congruence at $p$ as follows.
\begin{enumerate}
\item If $p \equiv 1 \mod 8$, then both $h_1$ and $h_3$ satisfy
the 3-term ASD congruence at $p$ given by $T^2 - sgn(p)a_1(p)T +
p^2$, where $sgn(p) = \pm 1 \equiv 2^{(p-1)/4} \mod p$~; \item If
$p \equiv 5 \mod 8$, then $h_1$ (resp. $h_3$) satisfies the 3-term
ASD-congruence at $p$
 given by $T^2 + 4ia_5(p)T - p^2$ (resp. $T^2 - 4ia_5(p)T - p^2$);
 \item If $p \equiv 3 \mod 8$, then $h_1\pm h_3$  satisfy the 3-term
 ASD congruence at $p$ given by
  $T^2 \pm 2\sqrt{-2}a_3(p)T -p^2$, respectively;
  \item If $p \equiv 7 \mod 8$, then $h_1\pm ih_3$  satisfy the
  3-term ASD congruence at $p$ given by
  $T^2 \mp 8\sqrt{-2}a_7(p)T -p^2$, respectively.
  \end{enumerate}
Here $a_1(p), a_3(p), a_5(p), a_7(p)$ are given by (17).
\end{theorem}

 Finally we remark that the character $\chi$ of
$G_{\Q(i)}$ may be viewed as an idele class character of $\Q(i)$ by
class field theory. Thus there is a cuspform $h(\chi)$ of $GL_2(\Q)$
whose local $L$-factors agree with those of $\chi$. The local
$L$-factors of $\rho \otimes \chi$ are in fact the local factors of
the form $f' \times h(\chi)$ on $GL_2(\Q) \times GL_2(\Q)$. We
summarize this discussion in

\begin{theorem} \label{GL(2)times GL(2)}
There are two cusp forms $f'$ and $h(\chi)$ of $GL_2(\Q)$
such that $\rho_-$ and $f'\times h(\chi)$ have the same local $L$-factors
over primes $p \ne l$.
\end{theorem}

Together with the fact that $L(s,\rho_+)$ is automorphic, we have shown

\begin{cor} There is an automorphic $L$-function over $\Q$ whose local
factors agree with those of the $l$-adic Scholl representation
 attached to the space $S_3(\G_4)$ at all primes $p
\ne l$.
\end{cor}

\section{Acknowledgements}
The authors are deeply indebted to the referee whose comments and
suggestions led to significant improvements of several proofs.
Special thanks are due to William A. Stein for facilitating part
of the computational results in this paper.

\end{document}